\documentclass[twocolumn]{autart}
\usepackage{amssymb,amsmath,dsfont}
\usepackage{graphicx}
\usepackage{subfigure,xcolor}
\usepackage{float}

\usepackage[english]{babel}
\usepackage{tikz,pgflibraryshapes,xcolor,xkeyval}
\usepackage{subfigure}
\usepackage{caption}

\newcommand {\zetavec}{{\boldsymbol{\zeta}}}

\newfont{\pseudocode}{cmtt10}

\newcommand{\bfx}{\mathbf{x}}
\newcommand{\bfl}{\boldsymbol{\ell}\,}

\newcommand{\bfB}{\mathbf{B}}

\newcommand{\bfA}{\mathbf{A}}

\newcommand{\bfC}{\mathbf{C}}

\newcommand{\bfg}{\mathbf{g}}

\newcommand{\bfG}{\mathbf{G}}

\newcommand{\dist}{\mathrm{dist}}

\newcommand{\pd}{{\partial}}

\newcommand{\Real}{\mathbb{R}}

\renewcommand{\Real}{\mathbb{R}}

\renewcommand{\Real}{\mathds{R}}

\newcommand{\beq}{\begin{equation}}
\newcommand{\eeq}{\end{equation}}

\newcommand{\bbm}{\begin{bmatrix}}
\newcommand{\ebm}{\end{bmatrix}}

\newcommand{\bpm}{\begin{pmatrix}}
\newcommand{\epm}{\end{pmatrix}}

\newcommand{\bit}{\begin{itemize}}
\newcommand{\eit}{\end{itemize}}

\newcommand{\ben}{\begin{enumerate}}
\newcommand{\een}{\end{enumerate}}

\newcommand{\barr}{\begin{array}}
\newcommand{\earr}{\end{array}}

\newtheorem{exam}{Example}[section]

\begin{document}
\begin{frontmatter}

\title{Supplementary material for:\\ Adaptive Observers and Parameter Estimation for a Class\\ of Systems Nonlinear in the Parameters}

\thanks[footnoteinfo]{This paper was not presented at any IFAC
meeting.  Cees van Leeuwen was supported by an Odysseus grant from the Flemish Science Organization FWO. Corresponding author I.~Yu.~Tyukin. Tel.
+44-116-2525106.}

\author[First]{Ivan Y. Tyukin}
\author[Second]{Erik Steur}
\author[Third]{Henk Nijmeijer}
\author[Second]{Cees van Leeuwen}

\address[First]{Dept. of Mathematics, University of Leicester, Leicester, LE1 7RH, UK (Tel: +44-116-252-5106; I.Tyukin@le.ac.uk)}
\address[Second]{Laboratory for Perceptual Dynamics, KU Leuven, Tiensestraat 102, 3000 Leuven, Belgium (erik.steur@ppw.kuleuven.be,  cees.vanleeuwen@ppw.kuleuven.be)}
\address[Third]{Dept. of Mechanical Engineering, Eindhoven University of Technology, P.O. Box 513 5600 MB,  Eindhoven, The Netherlands, (h.nijmeijer@tue.nl)}


\begin{keyword}
    Adaptive observers, nonlinear parametrization, weakly attracting sets
\end{keyword}

\begin{abstract}
    This supplement illustrates application of adaptive observer design from \cite{Tyukin:Automatica:2013} for systems which are not uniquely identifiable. It also provides an example of adaptive observer design for a magnetic bearings benchmark system \cite{ACC:2000:Lin}.
\end{abstract}

\end{frontmatter}

\section{Application of the design to non-identifiable systems}



\begin{exam}\label{example:non_ident:1}\normalfont Consider the following system
\begin{equation}\label{eq:exmaple:nonident:1}
\begin{array}{l}\dot{\bfx}=\bfA\bfx + \bfB \theta \lambda + \bfg(y),\\ y=\bfC^{T}\bfx\end{array},
\end{equation}
where
\[\ \bfA=\left(\begin{array}{cc} 0 & 1\\ 0& 0\end{array}\right), \ \bfB=\left(\begin{array}{c}1\\1\end{array}\right), \ \bfC=\left(\begin{array}{c}1\\0\end{array}\right), \ \bfg(y)=\left(\begin{array}{c}-2\\-1\end{array}\right) y
\]
and let $\Omega_\theta=[1,3]$, $\Omega_\lambda=[0.5,3]$,
$t_0=0$. In this particular case one can, in principle swap
$\theta$ and $\lambda$ (and $\Omega_\lambda$ with
$\Omega_\theta$ respectively).  It is clear that Assumption 3.1
is satisfied for this system and $y(t)$, as a function of $t$,
is bounded for all $t\geq t_0$. Moreover, the function $\varphi:$
$\varphi(\lambda)=\lambda$ and the function $\bfg$ are both continuous and
satisfy Assumption 3.2.

Let us now move to Assumption 4.1. Condition A1 is satisfied
for the given parametrization. Notice, however, that if the
domain $\Omega_\lambda$ contains $0$ then the condition would
not hold. In this situation if $0$ does not belong to the
domain $\Omega_\theta$ then swapping the definition of
$\lambda$ and $\theta$ (and $\Omega_\lambda$ with
$\Omega_\theta$) resolves the issue. Consider condition A2.
Sets $\mathcal{E}$ and $\mathcal{E}_0$, as follows from their
definition, coincide for this system (matrix $\bfG$ is
orthogonal to $\bfB$), and are defined as follows:
\[
\mathcal{E}_0(\lambda,\theta)=\mathcal{E}(\lambda,\theta)=\{(\lambda',\theta'), \lambda',\theta'\in\Real | \ \lambda\theta-\lambda'\theta'=0\}.
\]
We need to check if there is a function $\beta$  such that
\begin{equation}\label{eq:example:nonident:1:e_set}
|\lambda\theta-\lambda'\theta'| \geq \beta(\dist\left((\lambda',\theta'),\mathcal{E}(\lambda,\theta)\right)).
\end{equation}
In view of Remark 9, it is enough to show that for any
$(\lambda,\theta)\in\Omega_\lambda\times\Omega_\theta$ there is
a function $\beta$ (possibly dependent on $\lambda,\theta$)
such that (\ref{eq:example:nonident:1:e_set}) holds. And in
fact it is enough to show that it holds for all
$\lambda',\theta'$ from the domain
$\Omega_\lambda^\ast\times\Omega_\theta^\ast$ to which the
estimates $\hat{\lambda},\hat{\theta}$ may belong for $t\geq
t_0$. It is easy to see (from the proof of the theorem) that
Assumptions 3.1, 3.2, and A1 from Assumption 4. ensure that
both estimates $\hat{\lambda},\hat{\theta}$ are bounded for
$t\geq t_0$. Moreover, for the given set of initial conditions
for $\hat{\theta}$, one can estimate
$\Omega_\lambda^\ast\times\Omega_\theta^\ast$ a-priori. In this
example we set
$\Omega_\lambda^\ast\times\Omega_\theta^\ast=[0.5,3]\times[-10,10]$.

Let, for the sake of certainty, set $\lambda=2$, $\theta=1.5$.
In this case the set $\mathcal{E}$ consists of two isolated
curves. For any given point
$(\lambda_e,\theta_e)$ from
$\Omega_\lambda^\ast\times\Omega_\theta^\ast$  one can derive
two projections to the set $\mathcal{E}$. These projections
will be on the lines intersecting each branch orthogonally and
passing through $(\lambda_e,\theta_e)$. Let $d_e$ be the
minimal distance from $(\lambda_e,\theta_e)$ to these
projections. It is clear that for this distance $d_e$ there
will be a compact set $\mathcal{D}_e$ of points from
$\Omega_\lambda^\ast\times\Omega_\theta^\ast$ from which the
distance to $\mathcal{E}$ is exactly $d_e$. Given that
$\mathcal{D}_e$ is compact and that
$|\lambda\theta-\lambda_e\theta_e|$ is continuous one can
always define
\[
\vartheta: \ \vartheta(d_e)=\min_{(\lambda_e,\theta_e\in
\mathcal{D}_e
)}|\lambda\theta-\lambda_e\theta_e|.
\]
The function $\beta$ can now be defined as a function from
$\mathcal{K}_\infty$ satisfying: $\beta(d_e)\leq
\vartheta(d_e)$.

If finding the function $\vartheta$ analytically is complicated
then it can be estimated numerically. To illustrate the
procedure we randomly sampled
$\Omega_\lambda^\ast\times\Omega_\theta^\ast$ and derived the
values of $|\lambda\theta-\lambda_e\theta_e|$, $d_e$ for each
pair $(\lambda_e,\theta_e)$. The values of
$|\lambda\theta-\lambda_e\theta_e|$ then were plotted against
$d_e$ (please see the left panel in
\ref{fig:example:non_identifiabiliy:2}). The lower envelope of
the resulting set of points would then be an estimate of
$\vartheta$. One can observe that in this particular example
this envelope is already a strictly monotone function. This
clarifies condition A2.

\begin{figure}[!h]
\centering
\includegraphics[width=0.9\linewidth]{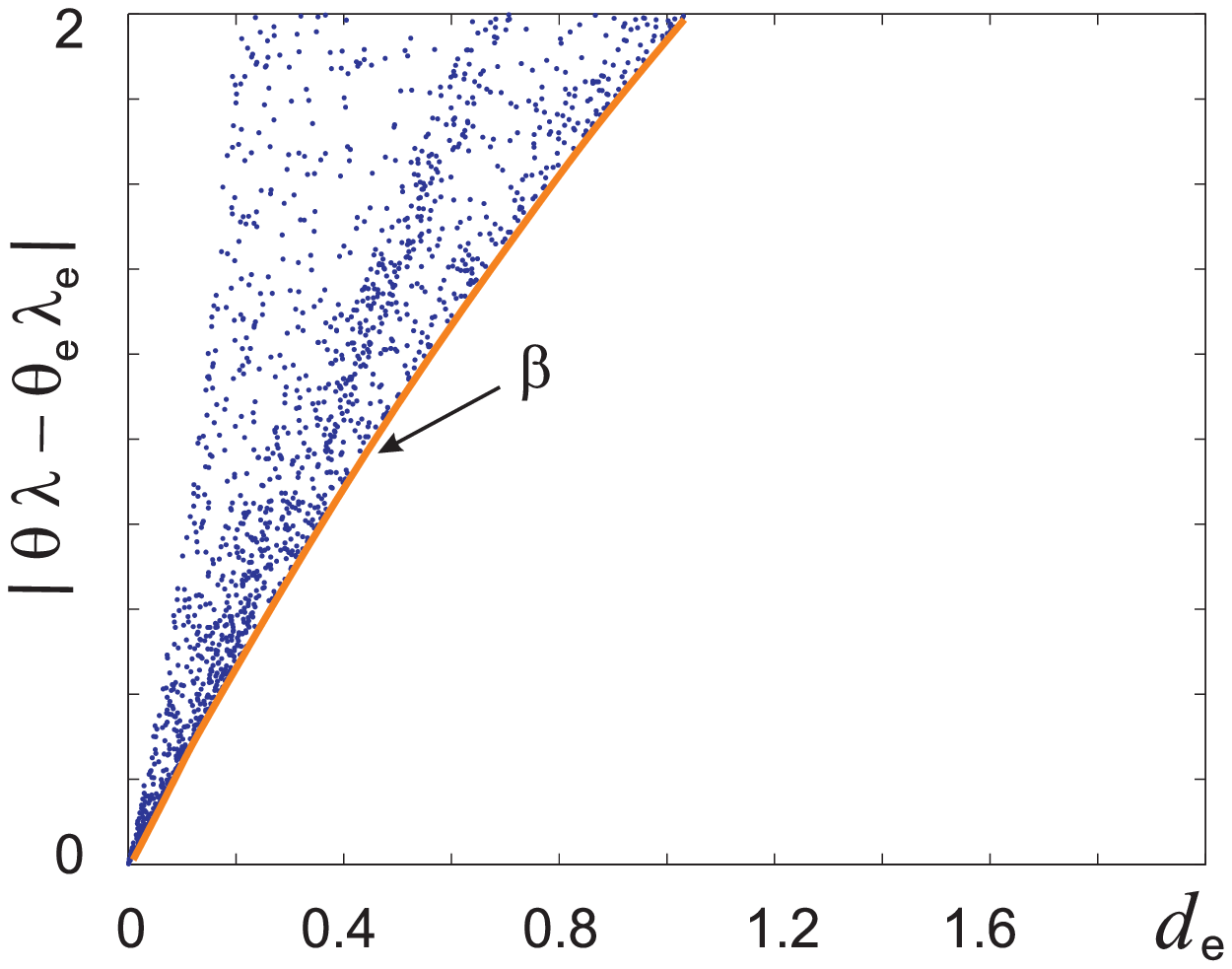}

\vspace{5mm}
\includegraphics[width=0.9\linewidth]{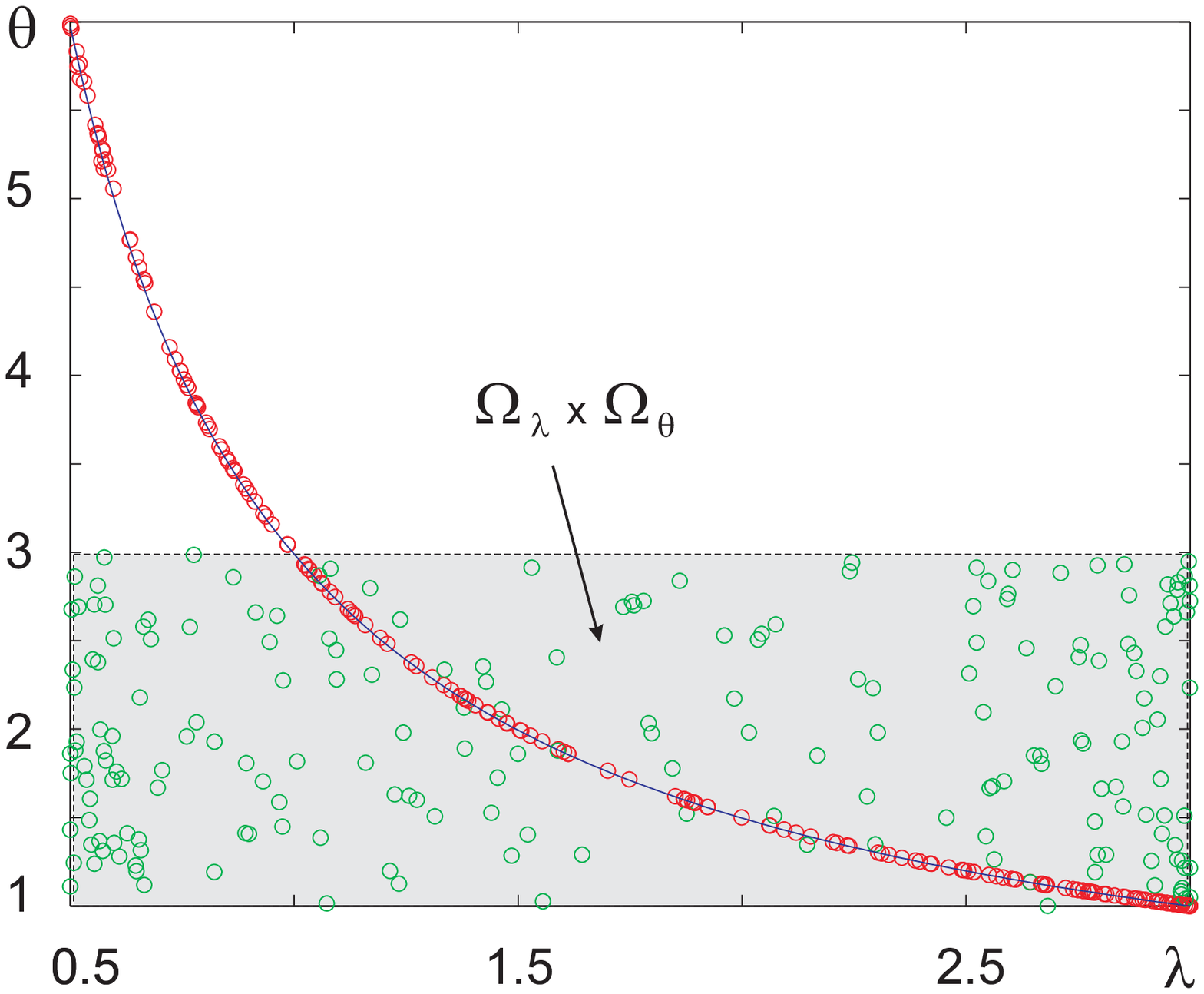}
\vspace{5mm}

\includegraphics[width=0.9\linewidth]{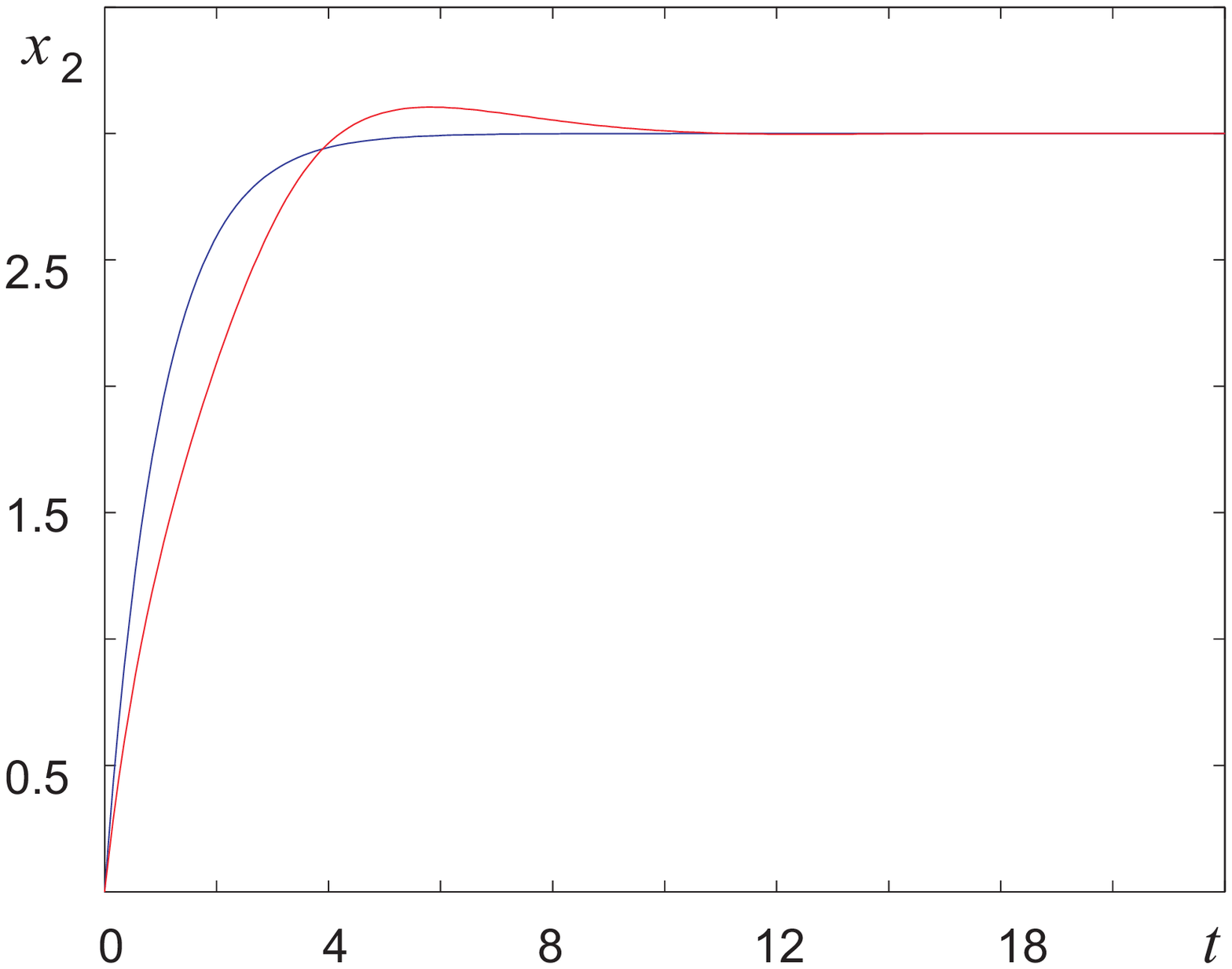}
\caption{Top panel: scatter plot of $10 000$ pairs ($|\lambda\theta-\lambda_e\theta_e|, d_e$) when $(\lambda_e,\theta_e)$ are sampled from the set $[-10,10]\times[0.5,3]$; red curve shows a boundary of this set  taken as an estimate of the function $\beta$. Central panel: red circles show the values to which the estimates $\hat\lambda,\hat{\theta}$ converge, green circles are the initial conditions, $\hat{\lambda}_0$, $\hat{\theta}_0$, and the corresponding branch of
 $\mathcal{E}$ is shown as the blue curve. Bottom panel: variable $x_2$ and $\hat{x}_2$ as functions of $t$.}\label{fig:example:non_identifiabiliy:2}
\end{figure}

Finally, we simulated the system together with the observer
(with $\gamma_\theta=1$, $\gamma=0.1$, $\bfl=(-2,-1)^{T}$) and
the results are summarized in Fig
\ref{fig:example:non_identifiabiliy:2}.

Notice that the vector of the estimates, $(\hat{\lambda},\hat\theta)$, may converge to points outside of
the domain $\Omega_\lambda\times\Omega_\theta$. It is, however, guaranteed that $\hat{\lambda}(t)\in\Omega_\lambda$ for all $t\geq t_0$.
Moreover, if
unmodeled dynamics is present, they may converge (depending on
initial conditions) to neighborhoods of elements from
$\mathcal{E}$ that are quite far from
$\Omega_\lambda\times\Omega_\theta$

\end{exam}

\begin{exam}\label{example:non_ident:2}\normalfont Consider now  system
\begin{equation}\label{eq:exmaple:nonident:2}
\begin{array}{l}\dot{\bfx}=\bfA\bfx + \bfB \theta + \bfg(y,\lambda),\\ y=\bfC^{T}\bfx\end{array} \  \bfg(y,\lambda)= \left(\begin{array}{c}0\\ e^{\lambda}\end{array}\right) + \left(\begin{array}{c}-2\\ -1\end{array}\right)y,
\end{equation}
where $\bfA$, $\bfB$, and $\bfC$, $\Omega_\theta$,
$\Omega_\lambda$, $\lambda$, $\theta$ are defined as in
(\ref{eq:exmaple:nonident:1}). It is clear that it satisfies
Assumptions 3.1, 3.2, and condition A1 of Assumption 4.1. Let
us check condition A2 of Assumption 4.1. First, we notice that
sets $\mathcal{E}$ and $\mathcal{E}_0$ are different in this
example. Indeed the set $\mathcal{E}_0$ consists of the single
element, $(\lambda,\theta)$, whereas the set $\mathcal{E}$ is
\[
\mathcal{E}(\lambda,\theta)=\{(\lambda',\theta'), \ \lambda',\theta'\in\Real | \theta + e^{\lambda}- \theta' - e^{\lambda'} =0 \ \}
\]
Similarly to the argument presented for the previous example,
one can conclude that the function $\beta$ (possibly dependent
on $\lambda,\theta$) exists, and can be estimated numerically
(see Fig. \ref{fig:example:non_identifiabiliy:3})
\begin{figure}[!h]
\centering
\includegraphics[width=0.9\linewidth]{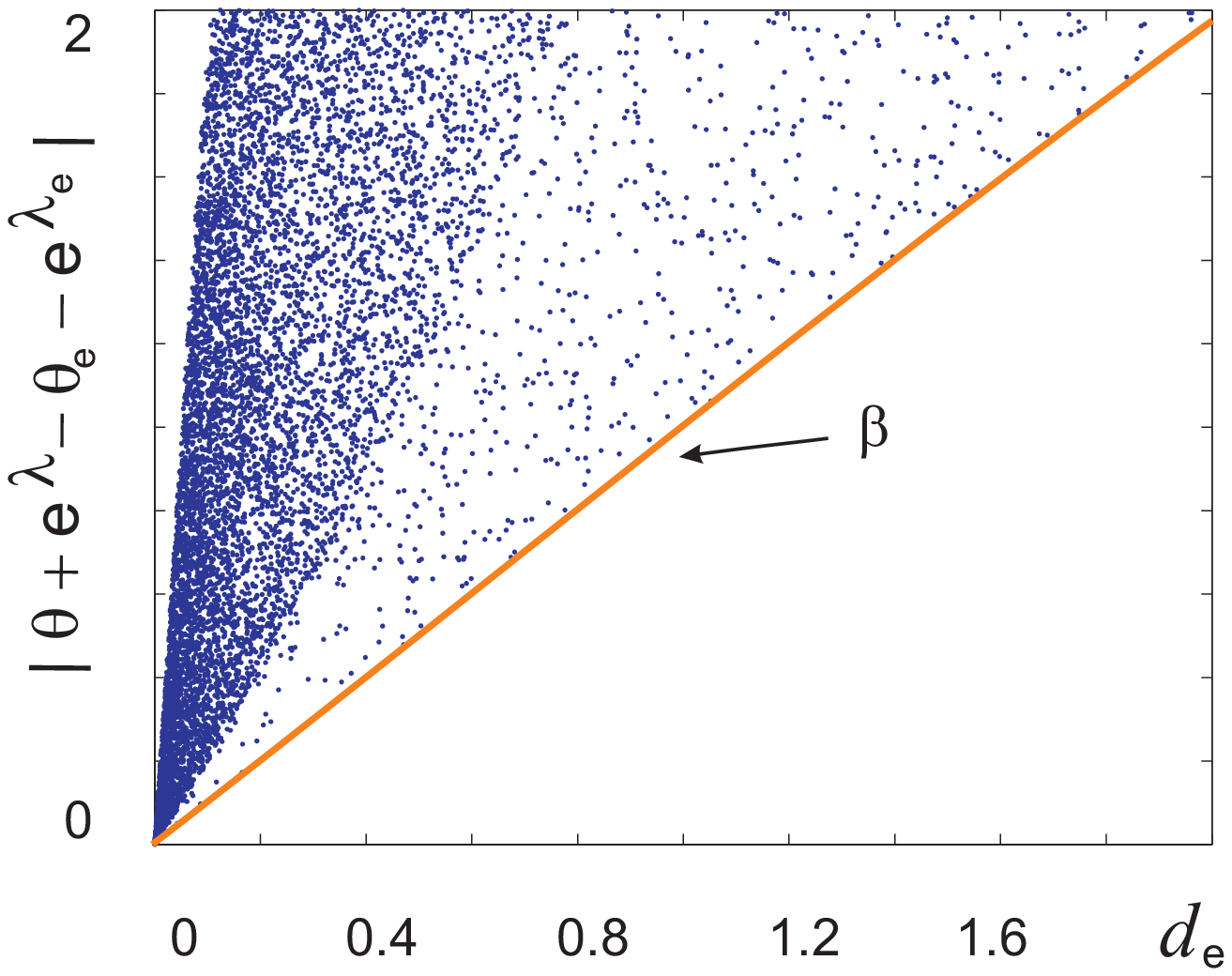}
 \vspace{5mm}
\includegraphics[width=0.9\linewidth]{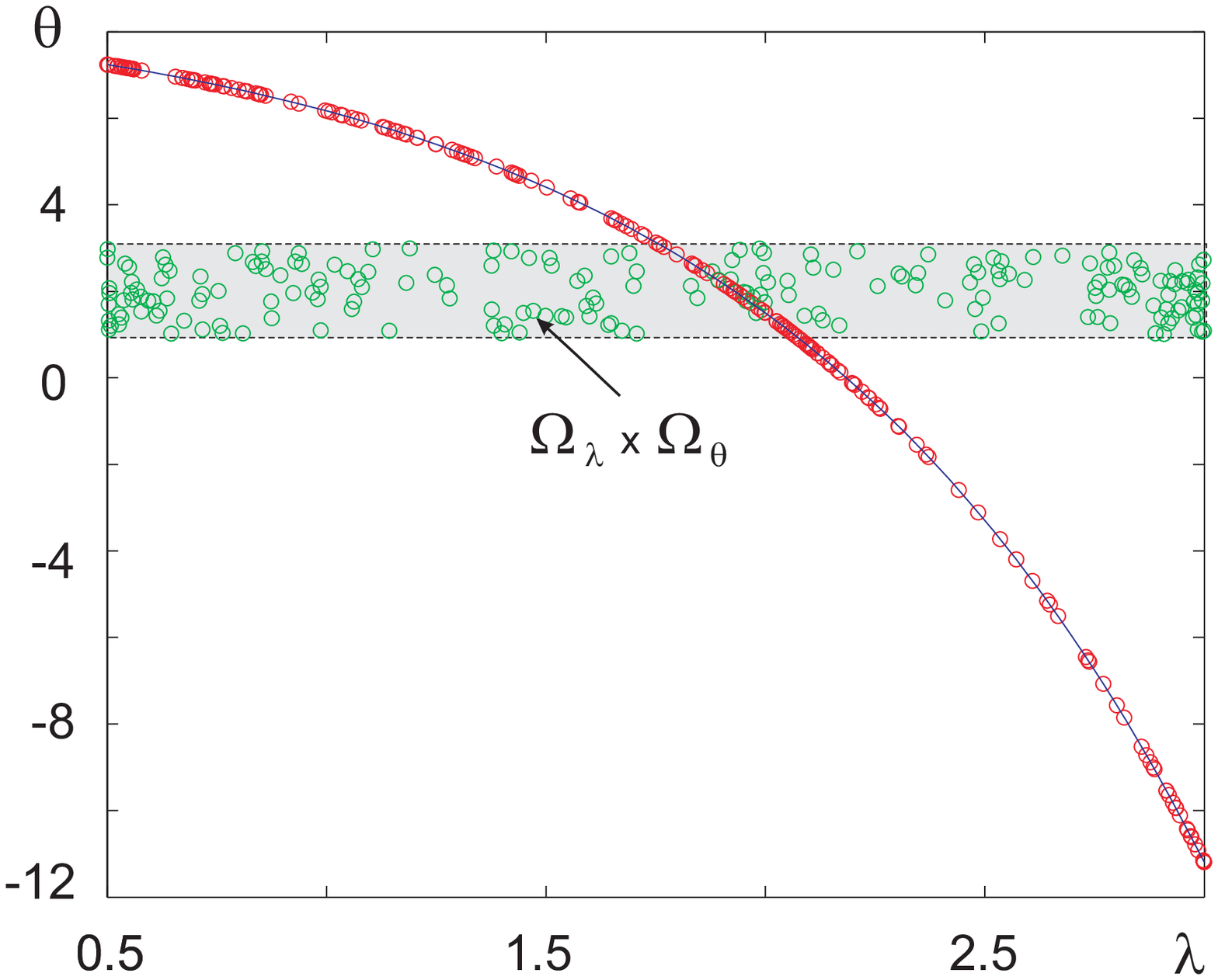}
\vspace{5mm}
\includegraphics[width=0.9\linewidth]{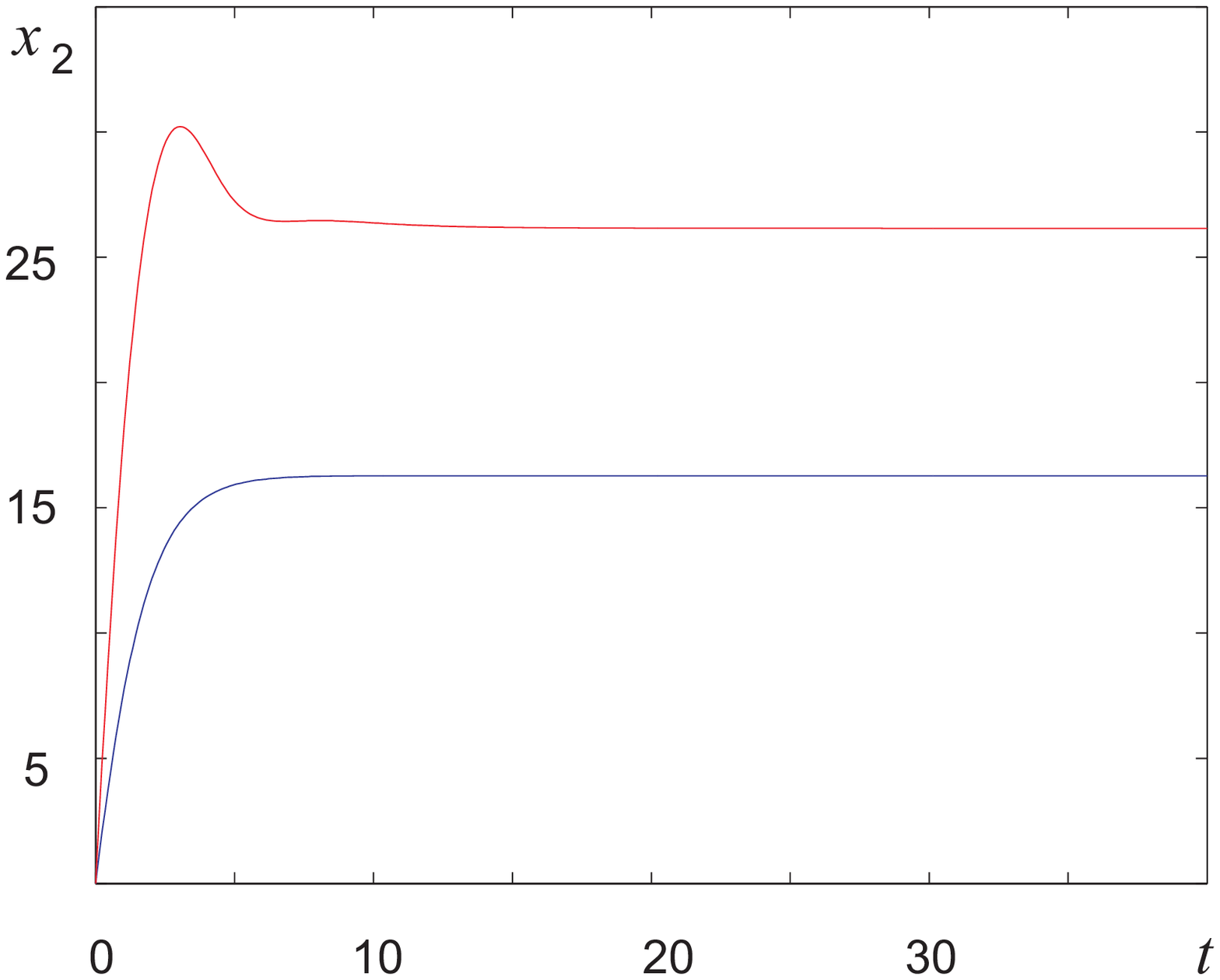}
\caption{Top panel: estimation of the function $\beta$ (obtained from $40 000$ pairs of points ($|\theta+e^{\lambda}-\theta_e-e^{\lambda_e}|, d_e$); $(\lambda_e,\theta_e)$ are sampled from the set $[-10,10]\times[0.5,3]$). Center panel:
 results of the estimation. Bottom panel: $x_2$ (blue) and $\hat{x}_2$ (red) as functions of $t$.}\label{fig:example:non_identifiabiliy:3}
\end{figure}
We simulated the system together with the observer
($\gamma_\theta=1$, $\gamma=0.05$, $\bfl=(-2,-1)^{T}$) for
various initial conditions in
$\Omega_\lambda\times\Omega_\theta$. The results are shown in
Fig. \ref{fig:example:non_identifiabiliy:3}. It is evident from
the figure that the state component $x_2$ is reconstructed with
an error, as expected (the Theorem guarantees reconstruction of
state only if $\mathcal{E}$ and $\mathcal{E}_0$ coincide).
\end{exam}

\section{State and parameter estimation for the magnetic beam benchmark system \cite{ACC:2000:Lin}}

The system dynamics is governed by the following equations \cite{ACC:2000:Lin} (see also Table 1, example 2):
\begin{equation}\label{eq:beam_example}
\dot{\bfx}=\left(\begin{array}{cc} 0 & 1\\ 0 & 0\end{array}\right) \bfx + \Psi(t,y(t),\lambda) \theta + \bfg(t), \ y=x_1,
\end{equation}
where
\[
\begin{array}{l}
\Psi(t,y(t),\lambda)=\left(\begin{array}{c} 0\\ 1\end{array}\right)\phi(t,y(t),\lambda); \ \bfg(t)=\left(\begin{array}{c}0\\ \frac{\sin(0.5t)}{1000}\end{array}\right)\\
\phi(t,y(t),\lambda)=\frac{cL}{J}  \left(\frac{q_2^2(t)}{(a(y(t)L+g)+\lambda b)^2}-\frac{q_1^2(t)}{(a(-y(t)L+g)+\lambda b)^2}\right),
\end{array}
\]
variables $q_1(t),q_2(t)$ are defined as
\[
\begin{array}{c}
\dot{q}_1=-\frac{R}{N} h q_1 + \frac{1}{N}\mathrm{sat}(u_1)\\
\dot{q}_2=-\frac{R}{N} h q_2 + \frac{1}{N}\mathrm{sat}(u_2)
\end{array}, \ \mathrm{sat}(u)=\mathrm{sign}(u)\min\{V_s,|u|\},
\]
and, $h$, the linear term of the otherwise nonlinear flux-current dependence, is modeled for simplicity as constant, $h=10^5$. Parameters $a,b,c,R,N,L,J,V_s$  are defined as in \cite{ACC:2000:Lin}, and stabilizing  control inputs $u_1, u_2$ are chosen as:
\begin{equation}\label{eq:beam_example_control}
\begin{array}{l}
2 x_1 + 2 \hat{x}_2 > 0 \Rightarrow \\
\ \  \left\{\begin{array}{l}u_1=0\\
u_2=N(\frac{R q_2}{N}-100 (q_2-\sqrt{\frac{2 x_1+2 \hat{x}_2J}{cL}+q_1^2}))\end{array}\right.
\\
2 x_1 + 2 \hat{x}_2 \leq 0 \Rightarrow\\
\ \ \left\{\begin{array}{l}u_2=0\\
u_1=N(\frac{R q_1}{N}-100 (q_1-\sqrt{\frac{2 x_1+2 \hat{x}_2J}{cL}+q_2^2}))\end{array}\right.
\end{array}
\end{equation}
\begin{equation}\label{eq:beam_example_high_gain_obs}
\left\{\begin{array}{ll}
\dot{\hat{x}}_1&=- 20(\hat{x}_1-x_1)+\hat{x}_2 \\
\dot{\hat{x}}_2&=-100(\hat{x}_1-x_1)+\phi(t,y(t),1)
\end{array}\right.
\end{equation}
The design is a combination of the high-gain observer (\ref{eq:beam_example_high_gain_obs}) with the switching controller (\ref{eq:beam_example_control}) aimed at ensuring that solutions of the unperturbed system obey linear second-order equation with asymptotically stable dynamics (cf. \cite{ACC:2000:Lin}).

Parameters $\theta$, $\lambda$ account for deviations of the true values of the magnetic constants $a,b,c$, and the air gap length, $g$, from their nominal values used in (\ref{eq:beam_example_control}), (\ref{eq:beam_example_high_gain_obs}). In particular, we suppose that the true values of $\theta,\lambda$ are unknown and belong to the interval $[0.8,1.2]$.

\begin{figure}
\centering
\includegraphics[width=0.9\linewidth]{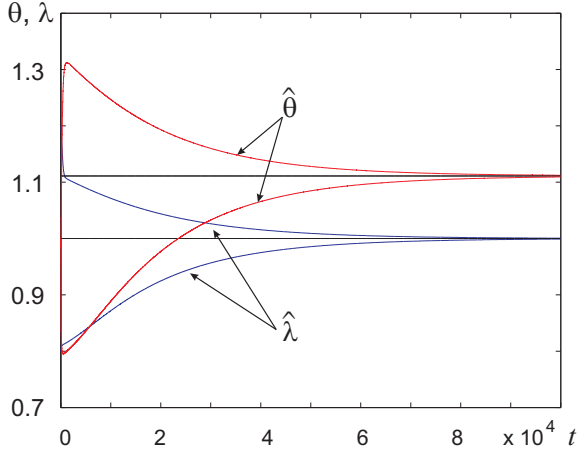}
\caption{Estimates $\hat{\theta},\hat{\lambda}$ of $\theta,\lambda$ and true values of $\theta,\lambda$ (black solid lines)}\label{fig:magnetic_beam_example}
\end{figure}

Let us now construct an observer for asymptotic reconstruction of $\theta,\lambda$. Equations (\ref{eq:beam_example}) are of the type (7), (71), and hence the observer candidate is defined by (72), (28)-(30). In particular, let $B=(1,1)^{T}$ and define $M(t,[\lambda,y])$, $\hat\zetavec$ as:
\begin{equation}\label{eq:beam_example_our_obs}
\begin{array}{l}
\dot{M}=\left(\begin{array}{cc}0 & 0\\ 0 & -1\end{array}\right) M + \left(\begin{array}{cc}0&0\\ -1 & 1\end{array}\right)\left(\begin{array}{c} 0 \\ \phi(t,y(t),\lambda)\end{array}\right)\\
\dot{\hat{\zetavec}}=\left(\begin{array}{cc}0 & 0\\ 0 & 1\end{array}\right)\hat{\zetavec}+\left(\begin{array}{c} -2\\-1\end{array}\right)((1,0)\hat\zetavec-y)\\
 \ \ \ \ \ \ + \ \bfB m_{21}(t,[\hat\lambda,y])\hat{\theta}+\bfg(t)\\
 \dot{\hat\theta}=-\gamma_\theta ((1,0)\hat\zetavec-y)m_{21}(t,[\hat\lambda,y])\\
 \hat{\lambda}=0.8+0.4(s_{1,1}+1)/2,
\end{array}
\end{equation}
where the variable $s_{1,1}$ is defined as in (79) with $\varepsilon=0$ (since no measurement noise and unmodeled dynamics are present). One can check that the vector
$(m_{21}(t,[\lambda',y]), R_1(t,\lambda',\lambda,\theta))$ (see the remark after Theorem 13),
where
\[
\begin{array}{l}
R_1(t,\lambda',\lambda,\theta)=\\
 \ \ \ \ \int_0^1 \int_{t_0}^t e^{-(t-\tau)} \frac{\pd }{\pd s}\phi(\tau,y(\tau),s(\xi,\lambda',\lambda))\theta d\tau d\xi, \\
 \ \ \ \  s(\xi,\lambda',\lambda)=\lambda' \xi + (1-\xi) \lambda,
\end{array}
\]
is $\lambda'$-UPE for $\lambda'\in[0.8,1.2]$, $\theta\in[0.8,1.2]$, and hence the system is uniquely identifiable, and assumptions of Theorem 13 hold. System (\ref{eq:beam_example}) together with observer (\ref{eq:beam_example_our_obs}) was simulated for various initial conditions with $\gamma_\theta=10^5$, $\gamma=150$, and samples of typical trajectories are shown in Fig. \ref{fig:magnetic_beam_example}.

\bibliographystyle{plain}
\bibliography{adaptive_observer_automatica}
\end{document}